\documentclass{article}
\usepackage{maa-monthly}

\theoremstyle{theorem}
\newtheorem{theorem}{Theorem}

\theoremstyle{definition}

\newtheorem*{remark}{Remark}
\newtheorem*{lemma}{Lemma}

\begin{document}

\title{An extension of Lobachevsky formula}
\author{Hassan Jolany}

\maketitle

\begin{abstract}
{In this paper we extend the Dirichlet integral formula of Lobachevsky. Let $f(x)$ be a continuous function and satisfy in the $\pi$-periodic assumption $f(x+\pi)=f(x)$, and $f(\pi-x)=f(x)$, $0\leq x<\infty $. If the integral $\int_0^\infty \frac{\sin^4x}{x^4}f(x)dx$ defined in the sense of the improper Riemann integral, then we show the following equality 
$$\int_0^\infty \frac{\sin^4x}{x^4}f(x)dx=\int_0^{\frac{\pi}{2} }f(t)dt-\frac{2}{3}\int_0^{\frac{\pi}{2} }\sin^2tf(t)dt$$ 

hence if we take $f(x)=1$, then we have $$\int_0^\infty \frac{\sin^4x}{x^4}dx=\frac{\pi}{3}$$ Moreover, we give a method for computing $\int_0^\infty \frac{\sin^{2n}x}{x^{2n}}f(x)dx$ for $n\in \mathbb N$}
\end{abstract}

\section{Introduction}

Dirichlet integral play an important role in distribution theory. We can see the Dirichlet integral in terms of distribution. The following classical Dirichlet integral has drawn lots of attention.

$$\int_0^\infty\frac{\sin x}{x}dx=\frac{\pi}{2}$$
We can use the theory of residues to evaluate this Dirichlet's integral formula.
G.H. Hardy and A. C. Dixon
gave a lot of different proofs for it. See [5-7]. In this paper we give an elegant method to generalize this Lobachevsky formula.

We start with the following elementary lemma. See [1-4]

\begin{lemma} For $\alpha\notin \mathbb Z\pi$, we have $$\frac{1}{\sin \alpha}=\frac{1}{\alpha}+\sum_{m=1}^{\infty}(-1)^{m}\left(\frac{1}{\alpha-m\pi}+\frac{1}{\alpha+m\pi}\right)$$
\end{lemma}
\begin{proof} For every positive integer $N$,
denote by $C_N$ the positively-oriented square in the complex plane with vertices $(N +\frac{1}{2})(\pm 1 \pm i)$. On the one hand, since the function $1/ \sin(z)$ is bounded on $C_N$
by a constant which is independent of $N$, one has $$\oint_{C_N}\frac{2\pi\alpha dz}{(z^2\pi^2-\alpha^2)\sin(\pi z)}\to 0$$ as $N\to \infty$. On the other hand, by the Residue Theorem, one also gets
$$\oint_{C_N}\frac{2\pi\alpha dz}{(z^2\pi^2-\alpha^2)\sin(\pi z)}=\sum_{n=-N}^N(-1)^n\frac{2\alpha}{n^2\pi^2-\alpha^2}+\frac{1}{\sin(\alpha)}$$

which proves the claim as $N \to \infty$.
\end{proof}

\begin{lemma}For $\alpha\notin \mathbb Z\pi$, we have the following identity 

$$\frac{1}{\sin^2 \alpha}=\frac{1}{\alpha^2}+\sum_{m=1}^{\infty}\left(\frac{1}{(\alpha-m\pi)^2}+\frac{1}{(\alpha+m\pi)^2}\right)$$
\end{lemma}
\begin{proof} The identity follows by differentiating termwise the classical
formula
, $$\cot(z)=\frac{1}{z}+\sum_{m=1}^{\infty}\frac{2z}{z^2-m^2\pi^2}$$
\end{proof}




\section{Lobachevsky formula} 

Now, we present the Lobachevsky formula.

\begin{theorem} Let $f$
satisfy the conditions of the beginning of the article. Then we have the following
Lobachevsky identity 

$$\int_{0}^\infty\frac{\sin^2 x}{x^2}f(x)dx=\int_{0}^\infty\frac{\sin x}{x}f(x)dx=\int_0^\frac{\pi}{2} f(x)dx$$

\end{theorem}

\begin{proof} Take $$I=\int_0^\infty \frac{\sin x}{x}f(x)dx$$

 we can write $I$ as follows 
$$I=\sum_{v=0}^{\infty}\int_{v\frac{\pi}{2}}^{(v+1)\frac{\pi}{2}}\frac{\sin x}{x}f(x)dx$$

where $v=2\mu-1$ or $v=2\mu$, by changing $x=\mu\pi+t$ or $x=\mu\pi-t$ we get
$$\int_{2\mu\frac{\pi}{2}}^{(2\mu+1)\frac{\pi}{2}}\frac{\sin x}{x}f(x)=(-1)^\mu\int_0^{\frac{\pi}{2}}\frac{\sin t}{\mu\pi+t}f(t)dt$$
and 

$$\int_{(2\mu-1)\frac{\pi}{2}}^{(2\mu)\frac{\pi}{2}}\frac{\sin x}{x}f(x)=(-1)^{\mu-1}\int_0^{\frac{\pi}{2}}\frac{\sin t}{\mu\pi-t}f(t)dt$$

so we get 

$$I=\int_0^{\frac{\pi}{2}}\frac{\sin t}{t}f(t)dt+\sum_{\mu=1}^{\infty}\int_0^{\frac{\pi}{2}}(-1)^\mu f(t)\left(\frac{1}{t+\mu\pi}+\frac{1}{t-\mu\pi}\right)\sin t dt$$

 Consequently we can write $I$ in the following form
 
 $$I=\int_0^{\frac{\pi}{2}}\sin t\left(\frac{1}{t}+\sum_{\mu=1}^{\infty}(-1)^\mu\left(\frac{1}{t+\mu\pi}+\frac{1}{t-\mu\pi}\right)\right)f(t)dt$$
 
 Hence $$I=\int_0^{\frac{\pi}{2}}f(t)dt$$ and proof of the identity $\int_{0}^\infty\frac{\sin x}{x}f(x)dx=\int_0^\frac{\pi}{2} f(x)dx$ is complete. Now we prove the second part of identity. Take $$J=\int_0^\infty \frac{\sin^2 x}{x^2}f(x)dx$$

 we can write $J$ as follows 
$$J=\sum_{v=0}^{\infty}\int_{v\frac{\pi}{2}}^{(v+1)\frac{\pi}{2}}\frac{\sin^2 x}{x^2}f(x)dx$$

where $v=2\mu-1$ or $v=2\mu$, by changing $x=\mu\pi+t$ or $x=\mu\pi-t$ we get
$$\int_{2\mu\frac{\pi}{2}}^{(2\mu+1)\frac{\pi}{2}}\frac{\sin^2 x}{x^2}f(x)=\int_0^{\frac{\pi}{2}}\frac{\sin^2 t}{(\mu\pi+t)^2}f(t)dt$$
and 

$$\int_{(2\mu-1)\frac{\pi}{2}}^{(2\mu)\frac{\pi}{2}}\frac{\sin ^2x}{x^2}f(x)=\int_0^{\frac{\pi}{2}}\frac{\sin^2 t}{(\mu\pi-t)^2}f(t)dt$$

so we get 

$$J=\int_0^{\frac{\pi}{2}}\frac{\sin^2 t}{t^2}f(t)dt+\sum_{\mu=1}^{\infty}\int_0^{\frac{\pi}{2}} f(t)\left(\frac{1}{(t+\mu\pi)^2}+\frac{1}{(t-\mu\pi)^2}\right)\sin^2 t dt$$

consequently we can write $J$ in the following form
 
 $$J=\int_0^{\frac{\pi}{2}}\sin^2 t\left(\frac{1}{t^2}+\sum_{\mu=1}^{\infty}\left(\frac{1}{(t+\mu\pi)^2}+\frac{1}{(t-\mu\pi)^2}\right)\right)f(t)dt$$
 
Hence from Lemma 1.2, we get $$J=\int_0^{\frac{\pi}{2}}f(t)dt$$ and proof is complete
\end{proof}

\section{Extension of the Lobachevsky formula}

Now we give a general method for calculating the following Dirichlet integral. 

$$\int_0^\infty \frac{\sin^{2n}x}{x^{2n}}f(x)dx$$

where $f(\pi+x)=f(x)$, and $f(\pi-x)=f(x)$, $0\leq x<\infty$. Here we have assumed $f$ is continuous and $\int_0^\infty \frac{\sin^{2n}x}{x^{2n}}f(x)dx$ defined in the sense of the improper Riemann integral. We start with $n=2$. As we did in the previous section, take $$I=\int_0^{\infty}\frac{\sin^4 x}{x^4}f(x)dx$$

By a direct computation $$\frac{d^2}{dx^2}\left(\frac{1}{\sin^2x}\right)=\frac{6}{\sin^4(x)}-\frac{4}{\sin^2(x)}$$

Next, differentiating twice termwise the right-hand side of identity of Lemma 1.2, we get the following identity 

$$\frac{1}{\sin^4\alpha}-\frac{2}{3\sin^2\alpha}=\frac{1}{\alpha^4}+\sum_{m=1}^\infty\left(\frac{1}{(\alpha-m\pi)^4}+\frac{1}{(\alpha+m\pi)^4}\right)$$

From the previous method which we explained in section 2, we can write $I$ as follows  

$$I=\int_0^{\frac{\pi}{2}}\sin^4t\left(\frac{1}{\sin^4t}-\frac{2}{3\sin^2 t}\right)f(t)dt$$

Hence 

$$I=\int_0^{\frac{\pi}{2}}f(t)dt-\frac{2}{3}\int_0^{\frac{\pi}{2}}\sin^2t f(t)dt$$

So we proved the following theorem 

\begin{theorem}  Let $f(x)$ satisfies in $f(x+\pi)=f(x)$, and $f(\pi-x)=f(x)$, $0\leq x<\infty $. If the following integral $$\int_0^\infty \frac{\sin^4x}{x^4}f(x)dx$$  defined in the sense of the improper Riemann integral, then we have the following equality 
$$\int_0^\infty \frac{\sin^4x}{x^4}f(x)dx=\int_0^{\frac{\pi}{2} }f(t)dt-\frac{2}{3}\int_0^{\frac{\pi}{2} }\sin^2tf(t)dt$$

\end{theorem}

As remark if we take $f(x)=1$, then we have 

\begin{remark} We have $$\int_0^\infty \frac{\sin^4x}{x^4}dx=\frac{\pi}{3}$$
\end{remark}

We have also the following remark from Lobachevsky formula

\begin{remark} If $f(x)$ satisfy in the condition $f(x+\pi)=f(x)$, and $f(\pi-x)=f(x)$, $0\leq x<\infty $, take $$I=\int_0^{\infty}\frac{\sin^{2n+1}x}{x}f(x)dx=\int_0^{\infty}\sin^{2n}x\frac{\sin x}{x}f(x)dx$$

If we set $\sin^{2n}xf(x)=g(x)$, we get $g(x+\pi)=g(x)$, $g(\pi-x)=g(x)$, now if we take $f(x)=1$, then $$\int_0^\infty \frac{\sin^{2n+1}x}{x}dx=\int_0^{\frac{\pi}{2}}\sin^{2n}xdx=\frac{(2n-1)!!}{(2n)!!}\frac{\pi}{2}$$

\end{remark}

Now, by the following important remark, we can calculate the Lobachevsky formula for any $n\geq 3$.  Let $f(z)$ satisfy the conditions of the beginning of the article.

\begin{remark} In fact, the Dirichlet integral $$\int_0^{\infty}\frac{\sin^{2n}z}{z^{2n}}f(z)dz$$
has the form (for $n\geq3$)
$$\alpha_1\int_0^{\frac{\pi }{2}}f(z)dz+\alpha_2\int_0^{\frac{\pi }{2}}\cot^{2n-2}(z)\sin^{2n-2}(z)f(z)dz+\cdots+\alpha_k\int_0^{\frac{\pi }{2}}\cot^{2}(z)\sin^{2}(z)f(z)dz$$

where the constants $\alpha_i$ can be computed by the use of the following formulas (and
the help of the engine Wolfram Alpha for instance) : For every positive integer $n$,
one can compute

$$\frac{d^n}{dz^n}\left(\frac{1}{\sin^2(z)}\right)=\frac{d^n}{dz^n}(1+\cot^2z)=\sum_{k=0}^n\binom{n}{k}\frac{d^{n-k}}{dz^{n-k}}(\cot z)\frac{d^k}{dz^k}(\cot z)$$

 by Leibnitz rule and then apply the closed formula $$\frac{d^m}{dz^m}(\cot z)=(2i)^m(\cot (z)-i)\sum_{j=1}^{m}{m \brace k}(i\cot(z)-1)^k$$ of Lemma 2.1. of [8], where ${n \brace k}$
are the Stirling numbers of the second kind. Now by applying the identity $$\frac{d^n}{dz^n}\left(\frac{1}{\sin^2z}\right)=\sum_{k=-\infty}^{\infty}\frac{(-1)^n(n+1)!}{(z+k\pi)^{n+2}}$$
we can find a closed formula for such Dirichlet integral formula for any $n$. For example, when $n=3$ we have $\{\alpha_1=\frac{2}{15}, \alpha_2=\frac{2}{15},\alpha_3=\frac{11}{15}\}$ and for $n=4$, $\{\alpha_1=\frac{272}{7!}, \alpha_2=\frac{64}{7!},\alpha_3=\frac{1824}{7!},\alpha_4=\frac{2880}{7!}\}$.
\end{remark}

\begin{remark}We have the following formulas

\;

1)$\int_0^{\infty}\frac{\sin^{6}z}{z^{6}}dz=\frac{11\pi}{40}$

2)$\int_0^{\infty}\frac{\sin^{8}z}{z^{8}}dz=\frac{151\pi}{630}$

3)$\int_0^{\infty}\frac{\sin^{10}z}{z^{10}}dz=\frac{15619\pi}{72576}$
\end{remark}

\section{Acknowledgment}I would like to thank of the anonymous reviewer for their careful reading of my manuscript and their many
insightful comments and suggestions and also correcting my explicit formula. 



%

\end{document}